\documentclass[submission,copyright,creativecommons]{eptcs}
\providecommand{\event}{ACT 2020} % Name of the event you are submitting to
\usepackage{underscore}           % Only needed if you use pdflatex.

\usepackage{amsmath}
\usepackage{comment}

\input diagxy

\newcommand{\cat}{\mathbf{cat}}
\newcommand{\mono}{\,\,\scalefactor{0.9}\to/ >->/\,\,}

\newcommand{\PF}{\vspace{\baselineskip}\noindent{\bf Proof. }}
\newcommand{\qed}{\hfill\rule{4pt}{8pt}\par\vspace{\baselineskip}}%Halbar

\newtheorem{definition}{Definition}

\newtheorem{proposition}[definition]{Proposition}

\newtheorem{remark}[definition]{Remark}

\def\mathopdef#1{\expandafter\def\csname#1\endcsname{\mathop{\rm#1}\nolimits}}
\def\mathoplsdef#1{\expandafter\def\csname#1\endcsname{\mathop{\rm#1}}}
\def\mathbfdef#1{\expandafter\def\csname#1\endcsname{{\rm\bf#1}}}
\def\mathrmdef#1{\expandafter\def\csname#1\endcsname{{\rm#1}}}

\newcommand{\bS}{\mathbf{S}}
\newcommand{\bT}{\mathbf{T}}

\newcommand{\bV}{\mathbf{V}}

\newcommand{\bX}{\mathsf{X}}
\newcommand{\bY}{\mathsf{Y}}
\newcommand{\bZ}{\mathsf{Z}}

\newcommand{\alens}{\mathsf{ALens}}
\title{The more legs the merrier: \\ A new composition for symmetric (multi-)lenses}

\author{Michael Johnson
\institute{CoACT, Departments of \\Mathematics and Computing\\
                Macquarie University }
\email{michael.johnson@mq.edu.au}
\and
Robert Rosebrugh
\institute{Department of Mathematics \\and Computer Science\\
                Mount Allison University}
\email{rrosebrugh@mta.ca}
}

\def\titlerunning{The more legs the merrier}
\def\authorrunning{Michael Johnson \& Robert Rosebrugh}

\begin{document}
\maketitle

\begin{abstract}
This paper develops a new composition of symmetric
lenses that preserves information which is important for implementing system
interoperation.
It includes a cut-down but realistic example of a multi-system business supply chain
and illustrates the new mathematical content with analysis of the systems, showing
how the new composition facilitates the engineering required to implement the
interoperations.
All of the concepts presented here are based on either pure category theory or
on experience in solving business problems using applied category theory.

\end{abstract}

\section{Introduction}
\label{sec-intro}
%(Note that
%we have 11 sections to fit into roughly 11 pages,
%so each one will be more succinct than we might like.)

Lenses are a category theoretic construct and are used in a very wide variety of applications.
Lenses come in a wide range of forms, but each kind of lens has a composition (associative, and with
identities), and so the
various lenses form the morphisms of categories, most often with objects which are themselves categories,
usually representing states
and transitions of some systems.  Among the kinds of lenses we will use here are asymmetric,
symmetric, and multiary lenses.

Symmetric lenses compose to, unsurprisingly, form new symmetric lenses.
Symmetric lenses are usually represented as spans of asymmetric lenses. Indeed, from the very
beginning, symmetric lenses have had various ad-hoc definitions, but in all cases the authors
noted that an alternative approach would be to define them as equivalence classes of spans of
asymmetric lenses.

In many applications, the fact that a symmetric lens might also be represented as
a \emph{co}span of asymmetric lenses is important, especially for implementation purposes.
However, the composition of symmetric lenses does not preserve the property that the lenses can be
represented by cospans --- two such symmetric lenses may (and frequently do) compose to form a
symmetric lens which cannot be represented as a cospan of asymmetric lenses.  Thus preserving the
factorisation to show how cospans of asymmetric lenses might be used in implementations becomes important.

In 2018, the first work on multilenses was begun.  Multilenses can be represented as multi-spans of asymmetric lenses
(often called wide spans, these multi-spans are spans with an arbitrary finite number of legs).
In this paper we analyse a small but realistic example of a supply chain in which the cospan representations
would be `composed away' by ordinary symmetric lens composition, and introduce a new kind of composition which
we call \emph{fusion} in which two ordinary symmetric lenses (spans with two legs) fuse to form a multilens with three
legs preserving the cospan representations, and more generally, two %symmetric
multilenses, spans with say $m$ and $n$ legs, fuse to form a %symmetric
multilens with $m+n-1$ legs, again preserving cospan representations.

% Hard coded section numbers for now!  Fix later.
The plan of the paper is as follows.  In Section 2 we present a cut-down example of actual supply chain
system interoperations.  Although the example is realistic, it has been cut to almost the very minimum
required to illustrate the mathematical developments in the rest of the paper.  The example is revisited
in several sections as we proceed through the development of the mathematics.  In Sections 3 and 4 we review
briefly asymmetric and symmetric lenses.  In the case of symmetric lenses the approach is the representation
of a symmetric lens as a span of asymmetric lenses.  In Section 5 we turn to cospans of asymmetric lenses,
pointing out the utility of a cospan representation and illustrating this with the example from Section 2.

Section 6 sketches some very new developments in lenses, the multilenses or wide spans of asymmetric lenses.
In Section 7 we introduce the new composition, \emph{fusion}, of multilenses, which in particular gives
a new composition for symmetric lenses.  We still refer to it as a kind of composition because it
contains all the information of the composite, but it has more, preserving, unlike normal compositions,
some important information about the makeup of the individual lenses that were composed.  This extra information
is shown in Section 8 to be just what is needed to preserve cospan representations when they exist, and so
to facilitate engineering practice.

After a brief interlude in Section 9 to review some other examples of fusion-like compositions, % and a
%place keeper for related work in Section 10,
we conclude and outline some future work in Section 10,
briefly describing how cospans among the feet of a multispan can be recovered unambiguously from
the appropriately fused multispan.  This emphasises the engineering importance of using fusion rather than
composition.

The paper requires only modest category theoretic background for which we refer readers to any standard
text.  One particular category theoretic notion we use repeatedly, above and below, is \emph{span}.  A \emph{span} in a category
is just a pair of arrows with common domain.  A \emph{cospan} is, of course, a pair of arrows with common codomain.
A \emph{wide span} is a collection of a finite number of arrows with common domain, and similarly a
\emph{wide cospan} is a collection of a finite number of arrows with common codomain.  We will frequently
talk about spans or cospans of asymmetric lenses.  Asymmetric lenses are defined in Section 3, but for now
note that an asymmetric lens is a functor, normally called the Get of the lens, and some further structure.
In our categories of lenses the arrows are lenses, and they are oriented in the direction of the lenses' Gets.
Thus a span of asymmetric lenses has, inter alia, two Gets (one for each lens) with a common domain.

%Remember to define span, cospan, point out that spans of lenses mean spans of Gs with lenses on them.  Similarly cospans.
%Include any other necessary category theoretic notations.

\section{An example}
\label{sec-eg}
%Describe ABC Frames, XYZ Warehouse, XYZ Sales, and XYZ Manufacturing.

We begin with an example which we will use to illustrate the concepts presented in this paper.  The example is based
on real systems, but they have been cut down to the essential details required to capture the ideas presented here.

Supply chains, especially global supply chains, are very much in the economic news at present because of the
disruptions to production and distribution caused by the corona virus crisis.  In many cases modern supply chains
are managed through system interoperations, with individual organisations owning and operating their own information systems,
but sharing enough information for those systems to interoperate automatically.  For example, two businesses, a supplier and a
customer, might both keep track of the amount of stock available on a customer's premises so that further stock
can be supplied in a just in time (JIT) manner.  To maintain the consistency of the two representations of the
customer's stock levels, information is exchanged between the two systems.  This is system interoperation in action,
and, as we will see later, it is mathematically captured by lenses --- certain kinds of bidirectional transformations.

\begin{figure}
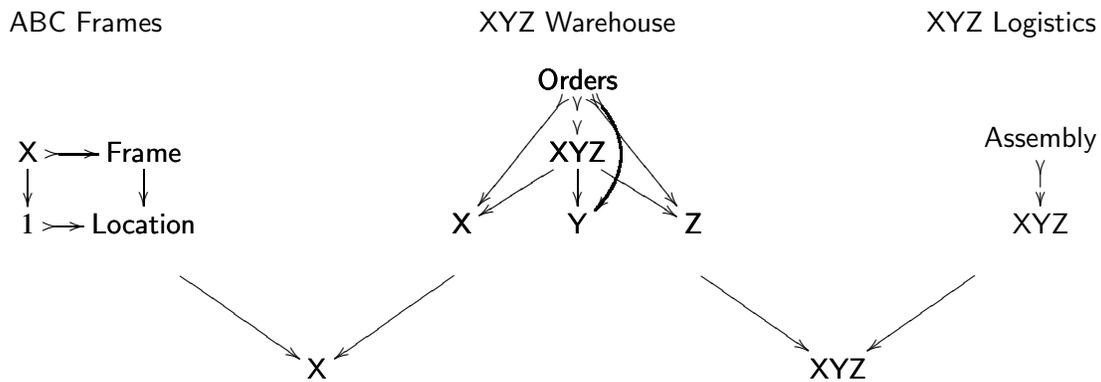

$$
\bfig
%% labels
\place(200,1200)[\mbox{\sf ABC Frames}]
\place(1900,1200)[\mbox{\sf XYZ Warehouse}]
\place(3400,1200)[\mbox{\sf XYZ Logistics}]
%% ABC Frames
\square(0,500)/ >->`->`->` >->/<400,250>[\bX`\mbox{\sf Frame}`1`\mbox{\sf Location};```]
%% Warehouse
\Atriangle(1500,500)/ >->` >->`/<400,500>[\mbox{\sf Orders}`\bX`\bZ;``]
\Atriangle(1500,500)/->`->`/<400,250>[\bX\bY\bZ`\bX`\bZ;``]
%\Vtriangle(0,0)/`->`->/<500,500>[\bA`\bC`\bV; `H'`G']
%\Vtriangle(1000,0)/`->`->/<500,500>[\bC`\bB`\bW; `K'`H'']
\morphism(1900,750)<0,-250>[\bX\bY\bZ`\bY;]
\morphism(1900,1000)/{@{ >-->}}/<0,-250>[\mbox{\sf Orders}`\bX\bY\bZ;]
\morphism(1900,1000)/{@{ >->}@/^1.5em/}/<0,-500>[\mbox{\sf Orders}`\bY;]
%%
%{\sf Assembly} \mono {\sf XYZ}
\morphism(3500,800)/{@{ >-->}}/<0,-300>[\mbox{\sf Assembly}`\bX\bY\bZ;]
%% left cospan
\morphism(1000,0)/{@{<-}}/<-500,350>[\bX`;]
\morphism(1000,0)/{@{<-}}/<500,350>[\bX`;]
%% right cospan
\morphism(2800,0)/{@{<-}}/<-500,350>[\bX\bY\bZ`;]
\morphism(2800,0)/{@{<-}}/<500,350>[\bX\bY\bZ`;]
\efig
$$
\caption{Three (simplified) business entities interoperating via two cospans of lenses}
\end{figure}

For our massively simplified businesses, imagine a supplier called ABC Frames.  It provides the basic structure, the
chassis if you will, for a variety of products that are manufactured by other businesses which make up ABC Frames'
customers.  One of those customers is XYZ Manufacturing, which is operated as several distinct business entities,
XYZ Warehouse, XYZ Logistics and XYZ Production being the businesses we will consider here.

A true representation of ABC Frames will, like most businesses, involve systems which track hundreds, and often many
thousands of different types of entities, along in each case with their many attributes (for example, colour, dimensions,
location, serial number, base price, etc).  For our purposes we focus on a single type of frame, called here just
{\sf Frame}s.  These frames almost certainly have many attributes including those above, but for now we will just consider location.
(Other attributes can be easily managed too, but will just clog up our pictures if we depict them.)
So, for our purposes we might think of ABC Frames as having at any instant a set of frames, with each frame having a specified
location.  For our purposes ABC Frames' information system stores two sets and a function between them
$$
{\sf Frame} \longrightarrow {\sf Location}.
$$

XYZ Manufacturing, in its simplified form, takes frames which they call {\sf X}, and assembles them with other
products called simply {\sf Y} and {\sf Z} to produce a consumer product known as XYZ.  In particular, XYZ Warehouse
keeps stocks of {\sf X}, {\sf Y} and {\sf Z}, along with various attributes of each instance of those stocks (their serial
numbers, colours, and so on) which we will not record here.  The warehouse information system also calculates from
its known stocks the product ${\sf X} \times {\sf Y} \times {\sf Z}$ which is useful for XYZ's consumer facing operations
because it shows all of the possible available combinations of {\sf X}, {\sf Y} and {\sf Z} items that might be
assembled and provided to a potential customer.  XYZ Warehouse's information system also stores briefly some
information about orders, but we will return and fill that in later.  So, for now,
mathematically speaking, XYZ Warehouse's information systems stores
three sets, along with their product and the product projections (all of which can be calculated as required from
the three sets).

XYZ Logistics, also known as XYZ Sales because it is the principal interface with XYZ's customers, assembles orders
placed by customers after, of course, detailed discussions with those customers about customers' needs and desires,
frequently consulting XYZ's current catalogue, which is really just a copy called {\sf XYZ} of the product
${\sf X} \times {\sf Y} \times {\sf Z}$ calculated by XYZ Warehouse.
So XYZ Logistics' information system stores a set whose value is kept consistent with
XYZ Warehouse's calculated product by system interoperations (a lens in fact).
In addition, XYZ Logistics keeps track of customer orders, so it maintains a set usually called {\sf Assembly Order}
(or just {\sf Assembly} for short)
whose elements are casually called order-lines or order items, which should be a subset of {\sf XYZ} (a subset
because, of course, we don't want two orders for the same assembled product --- each assembly is unique and can
only be sold once).  As you can imagine, XYZ Logistics does much more, and has its ``Logistics'' name, and uses terms like
{\sf Assembly} for orders, because it looks after many logistical issues including the transfer of products
from the warehouse to XYZ Production, but we don't need to discuss these things here and they have been elided.
Mathematically, the much simplified information system for XYZ Logistics contains a set {\sf XYZ} and a subset
of that set representing the current order items:
$$
{\sf Assembly} \mono {\sf XYZ}.
$$

The eagle-eyed among readers will have spotted another concern.  Since the catalogue, {\sf XYZ} contains all
possible assemblies from products contained in the warehouse, customers might order distinct assemblies which
nevertheless contain the same instance of a particular product.  For example, two different customers might
order two different assemblies both of which are built on the frame with serial number 4097.  We need to guard
against this as both orders can't be satisfied --- there is only one frame 4097.   This is where the ``extra
information'' mentioned above in discussing XYZ Warehouse's information system comes into play.  XYZ Warehouse
also keeps a local copy called {\sf Orders} of XYZ Logistics' {\sf Assembly}, maintained again by system interoperations, along with
the corresponding subset inclusion of {\sf Orders} into ${\sf X} \times {\sf Y} \times {\sf Z}$.  Since that product
comes with its product projections in XYZ Warehouse, the information system there can see the composite of the inclusion
with each of the product projections, and, as part of its inbuilt constraints, it requires that those compositions are
monic (injections).  Any attempt to enter an order item that violates that constraint will be rejected, and
via the systems interoperations, a customer will be unable to order an assembly containing, for example, frame 4097
if there is already another extant order for an assembly using that frame.  Mathematically the full version of
the fragment of XYZ Warehouse's information system that we will be considering is summarised in the diagram
below.
$$
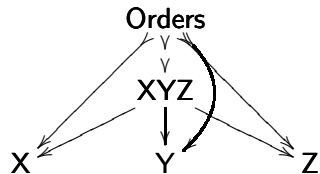
\bfig
\Atriangle(0,0)/ >->` >->`/<500,500>[\mbox{\sf Orders}`\bX`\bZ;``]
\Atriangle(0,0)/->`->`/<500,250>[\bX\bY\bZ`\bX`\bZ;``]
%\Vtriangle(0,0)/`->`->/<500,500>[\bA`\bC`\bV; `H'`G']
%\Vtriangle(1000,0)/`->`->/<500,500>[\bC`\bB`\bW; `K'`H'']
\morphism(500,250)<0,-250>[\bX\bY\bZ`\bY;]
\morphism(500,500)/{@{ >-->}}/<0,-250>[\mbox{\sf Orders}`\bX\bY\bZ;]
\morphism(500,500)/{@{ >->}@/^1.7em/}/<0,-500>[\mbox{\sf Orders}`\bY;]
\efig
$$

That completes our summary of XYZ Manufacturing's information systems, and we have seen the very simple information
systems maintained by ABC Frames, but to complete our automated supply chain we need to see how the two companies'
systems interact, and it is very simple.  Among the locations where ABC Frames might keep track of frames is
XYZ Warehouse.  So in ABC Frames there is an element $1 \longrightarrow {\sf Location}$ which picks out XYZ's
warehouse, and the pullback
$$
\bfig
\square(0,0)/ >->`->`->` >->/<400,250>[\bX`\mbox{\sf Frame}`1`\mbox{\sf Location};```]
\efig
$$
calculates the subset of {\sf Frame}s which are the frames located at XYZ Warehouse.  That should of course
correspond to {\sf X} in XYZ's own systems and system interoperations are used to keep those two sets consistent.
This supports for these two companies their version of the JIT supply system described at the beginning of this section.

Although this example is vastly simplified, it does model many interesting aspects of category theoretic
information systems interoperation, including, as we will see below, symmetric lenses, multilenses, cospan
implementations of interoperations, amendment lenses, and so on.  A summary of the three simplified business
entities alonng with rough indications of the lenses between them used for maintaining interoperations (excluding
part of the amendment lens synchronising {\sf Asssembly} and {\sf Orders}) is shown in Figure 1.

Before saying more about all this
we review the relevant concepts from earlier work, and develop the new theory required for this paper.

\section{Asymmetric lenses}
\label{sec-asymlens}
%Clarke definition.  Rough connection to earlier definitions.  Sufficiency via our unification paper.
%Relationship to set-based lenses and Haskell lenses could be done here, but isn't at this stage.

Lenses are used to maintain synchronisation between (or in the case multilenses, among) different
systems.  In asymmetric lenses, one of the systems (the one with state space $S$ below) has all the information
required to reconstruct the other (the one with state space $V$ below).  An operation, usually called ``Get'' and frequently
denoted $G$, gives for any state of the system $S$, the corresponding state of the system $V$.
In the reverse direction, we would not expect a $V$ state to contain enough information to recreate an entire $S$ state.
Instead, the operation usually called ``Put'' provides a new state $s'$ of $S$ given an old state $s$ of $S$ and a change of state
in $V$ from $Gs$ to some new state $v'$, such that the new states of $S$ and $V$, $s'$ and $v'$, are again synchronised.

Naturally, state spaces will be represented here as categories --- a state of a system is an object of the
state space, and arrows of the state space are state transitions (and state transitions can be composed associatively,
and there are identity transitions corresponding to no-change).  Thus the state transition in $V$ just mentioned is an
arrow $Gs \to v'$ in $V$.

Database view updating \cite{bsusrv} provides a typical (and longstanding) example:  Suppose that $S$
is the state space of the information system of ABC Frames, one of the organisations discussed in the previous
section.   The object $\bX$ of the previous section is the state space of a view of $S$.  A state of $\bX$ is
just a set, the current set of frames located at the XYZ Warehouse according to the current state of the information
system $S$.   And the pullback above shows how to calculate $\bX$ from a current state of ABC Frame's information system
(as discussed in Section~\ref{sec-eg}, a state of ABC Frame's information system is just a function
${\sf Frame} \longrightarrow {\sf Location}$).  The Get of this view is calculated by the pullback, which in
database terms is simply the query ``{\sf select Frame where Location equals XYZ Warehouse}''.

If the view $\bX$ is changed, then the Put needs to construct sets {\sf Frame} and {\sf Location} and a function
${\sf Frame} \longrightarrow {\sf Location}$.  The most natural choice of Put in this case starts from the old
function ${\sf Frame} \longrightarrow {\sf Location}$, leaves {\sf Location} unchanged, changes {\sf Frame} to
correspond to the new $\bX$ by adding or deleting elements as required,
retains the values of the function for all those elements of Frame that remain in
the new $\bX$, and assigns any extra elements of $\bX$ to have location XYZ Warehouse (since if the set $\bX$
is intended to be the result of the above query, any extra (new) elements in $\bX$ can be assumed to located
at XYZ Warehouse).

All of this can be formalised easily using the theory of database modelling via EA-sketches \cite{jrfuvup} %Check citation
in which diagrams like those from the previous section are the base graphs of sketches, limits are used to
ensure that things like products and monics are appropriately realised, colimits are used to define attributes,
and the state spaces just described are models of the sketches, that is, full subcategories of finite set-valued
functor categories for which the functors preserve finite limits and finite coproducts.

We turn now to the
formal definition of asymmetric lenses.  For readers who are most familiar with early work on lenses \cite{pslvut, olesthes}
or with lenses as implemented in Haskell, this definition might come as a surprise, but it elegantly
captures the generality required in clear category theoretic terms.  Asymmetric lenses, as defined here, are sometimes
called d-lenses or delta lenses \cite{ddl}, and unify a wide range of different types of lenses
\cite{jrusbdbebl}.

\begin{definition}
{\rm (Clarke \cite{BryceACT19}):}  An \emph{asymmetric lens} is a commutative triangle of functors, as depicted below left,
in which $F$ is a discrete opfibration, $P$ is bijective on objects,
and $G$ and $P$ are called the Get and the Put respectively.
%Half-Bowtie (= triangle) and composable lenses as concatenated triangles
$$
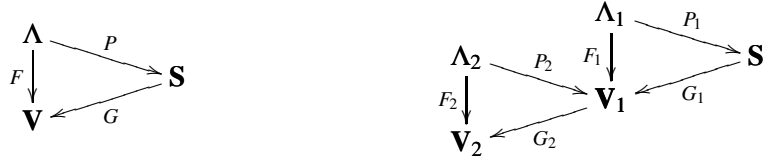
\bfig
\Dtriangle(0,75)/->`->`->/<500,150>[\Lambda`\bS`\bV;F`P`G]
%\Ctriangle(500,0)/->`->`->/<500,150>[\bD`\bC`\bE;``]
\Dtriangle(1500,0)/->`->`->/<500,150>[\Lambda_2`\bV_1`\bV_2;F_2`P_2`G_2]
\Dtriangle(2000,150)/->`->`->/<500,150>[\Lambda_1`\bS`\bV_1;F_1`P_1`G_1]
\efig
$$
\end{definition}

The category $\Lambda$ is, up to equivalence, a category
with the same objects as $\bS$ and with an arrow $\alpha$ from $s$ to $s'$ if and only if
$\alpha : Gs \to Gs'$ in $\bV$ (using the notation $s'$ from the beginning of this section).
The span $(\Lambda,F,P)$ is a co-functor \cite{Aguiar} from $\bV$ to $\bS$.
For further motivation and details we refer the reader to \cite{BryceACT19}.
As noted there, composition of asymmetric lenses is defined by simply composing the Gets and pulling back the
second Put, $P_2$, along the first lens's discrete opfibration, $F_1$ (see the above right diagram).
Thus, there is a category $\alens$ whose objects are categories and whose arrows are asymmetric lenses, oriented in the
direction of the Get (so the above left triangle is an arrow of $\alens$ from $\bS$ to $\bV$).

\section{Symmetric lenses as spans of asymmetric lenses}
\label{sec-symlens}
%Mention equivalence, but suppress it for this paper.
%Show how LR behaviour comes from Clarke definition.  Cite Clarke ACT submission.
%Composition of symmetric lenses, either by ``pullback'', or by composing such LR behaviours.
%Example between ABC Frames and XYZ Warehouse will be saved for the next section.

Lenses, as just defined, are examples of bidirectional transformations \cite{OxfordNotes}.  To
reiterate, a bidirectional transformation maintains consistency between two systems as one or the
other changes, and the functor part and the cofunctor part of an asymmetric lens embody the two updates required,
one in each direction, to restore consistency after a change of state of one system or the other.
As we've noted, such lenses are often
called \emph{asymmetric} lenses to emphasise the asymmetry noted at the beginning of Section~\ref{sec-asymlens}:
A state of one system, $\bS$, has all the information required to construct a state of the other system,
$\bV$, and this is reflected in the fact that one of the updates, $G$, is simply a functor.

While asymmetric lenses do arise in real world applications of bidirectional transformations,
there are many important cases where neither system has the information to reconstruct the other
completely.  Instead, each system ``knows'' things that the other system does not.  What's required
is a \emph{symmetric} lens \cite{hpwsl,dsl}.  As was conjectured in both the papers just cited, and
in \cite{dmaccat12},
a symmetric lens can be defined as an equivalence class of spans of asymmetric lenses \cite{jrjot}.
In this paper we will elide the details about the equivalence
(full details are available in \cite{jrjot})
and work with representatives of equivalence classes.
Again, the best available modern treatment is due to Bryce Clarke.

\begin{definition}
\label{def-symlens}
{\rm This formulation is due to Clarke \cite{BryceACT20}:}
A (representative for a) \emph{symmetric lens} is a span of asymmetric lenses as shown,
%Bowtie
$$
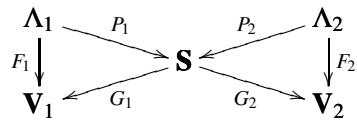
\bfig
\Dtriangle(0,0)/->`->`->/<500,150>[\Lambda_1`\bS`\bV_1;F_1`P_1`G_1]
\Ctriangle(500,0)/->`->`->/<500,150>[\Lambda_2`\bS`\bV_2;P_2`F_2`G_2]
\efig
$$
in which the objects are categories, the arrows are functors, the vertical arrows $F_1$ and $F_2$ are discrete opfibrations,
and the functors $P_1$ and $P_2$ are bijective on objects.
\end{definition}

This ``bowtie'' representation of symmetric lenses turns out to be particularly convenient.  For example,
as a bidirectional transformation, a symmetric lens should show how to restore consistency if a state of
either $\bV_1$ or $\bV_2$ is changed.  These two operations have variously been called the Rightward and Leftwards
\cite{hpwsl} and Forwards and Backwards \cite{dsl} propagations.  Each propagation is easily visible in the
bowtie, with, for example, the Forwards propagation given by the span $(\Lambda_1,F_1,G_2 P_1)$, or in short,
the South-East diagonal $G_2 P_1$.  In more detail:  The systems with state spaces $\bV_1$ and $\bV_2$ are synchronised
when there is an $s$ in $\bS$ with the current states of $\bV_1$ and $\bV_2$ equal to $G_1 s$ and $G_2 s$ respectively.
If $\bV_1$ then changes state via say an $\alpha : G_1 s \to v_1 '$, then that determines a unique arrow of $\Lambda_1$,
$\hat{\alpha}: s \to s'$ and $G_2 P_1 \hat{\alpha}$ is an arrow of $\bV_2$, of the form $G_2 s \to v_2 '$.
Furthermore, the new states $v_1 '$ and $v_2 '$ are synchronised by $s'$. Synchronisation has been restored by the
Forward propagation $G_2 P_1$.

This might be a convenient moment to say a few more words about the equivalence relation that we are
mostly suppressing in this short version to ease the reader's burden.  A symmetric lens is a bidirectional
transformation between $\bV_1$ and $\bV_2$, while $\bS$ is generally considered to be hidden coordination information.
If two bowtie representations between $\bV_1$ and $\bV_2$ have the  same Forwards and Backwards propagations
then they should be considered to be representatives of the same abstract symmetric lens even if they happen
to manage their coordination via different categories $\bS_1$ and $\bS_2$.  Again, we refer readers interested
in full details about the required equivalence relation to \cite{jrjot}. %Check!

Another convenience of the bowtie representation is that it shows immediately how symmetric lenses compose.
In fact, it shows that in two distinct, but equivalent ways.  First, operationally, the Forward propagation
described two paragraphs ago results in an arrow of $\bV_2$, which can be in exactly the same manner Forward
propagated along a symmetric lens from $\bV_2$ to $\bV_3$ defining a composite Forward propagation.  Similarly,
Backward propagations can be iterated thus determining, operationally, a composite symmetric lens from $\bV_1$ to
$\bV_3$.  This can be shown to be equivalent to the composite span of asymmetric lenses described in Remark~\ref{rem-spancompn}
below using lens structures on the pullback in $\cat$.  First we describe these ``pulled back'' lenses and their
basic properties.

\begin{proposition}
\label{prop-lenspb}
Given a cospan of asymmetric lenses as shown
% Reverse-Bowtie
$$
\bfig
\Dtriangle(800,0)/->`->`->/<500,150>[\Lambda'_1`\bS'`\bV'_1;F'_1`P'_1`G'_1]
\place(650,0)[=]
\Ctriangle(0,0)/->`->`->/<500,150>[\Lambda_2`\bS`\bV_2;P_2`F_2`G_2]
\efig
$$

\begin{enumerate}
\item Each lens pulls back along the other lens's Get to give a lens
\item The resulting square of asymmetric lenses commutes in $\alens$, and will be referred to as the ``pullback'' of the cospan
\item The cospan itself determines operationally Forwards and Backwards propagations
\item And the propagations determined by the cospan and by its ``pulled back'' span coincide.
\end{enumerate}
\end{proposition}

\PF
We give brief proof outlines:
\begin{enumerate}
\item This is proved by explicitly constructing the Puts in \cite{jrjot}, but here we just outline the simple proof due
to Clarke \cite{BryceACT19}: %Check!
Part~1 follows immediately from the pullback pasting lemma and the facts that
discrete opfibrations pull back along functors to give discrete opfibrations
and bijective on objects functors pullback along functors to give bijective on object functors.
\item The Gets of the two sides of the square commute by construction (the Gets form a pullback square in $\cat$)
and it's easy to see from the explicit construction of the Puts along the two sides of the square that
the compositions of the cofunctors coincide too.
% Is there an easy cofunctor proof? (avoiding our nuts-and-bolts construction)
\item For a cospan labelled as above, we call objects $s$ of $\bS$ and $s'$ of $\bS'$ \emph{synchronised} when $G_2 s = G'_1 s'$.
If $s$ and $s'$ are synchronised, and $\alpha : s \to r$ is a change of state of the system with state space $\bS$, then
the Forward propagation of $\alpha$ is the Put, $P'_1$, of $G_2 \alpha : G'_1 s' \to G_2 r$,
which will be an arrow of $\bS'$ with domain $s'$ and whose codomain is then synchronised with $r$.
The Backward propagation is defined similarly.
\item  Finally, by the constructions of pullbacks in $\cat$, $s$ and $s'$ are synchronised by the pulled back span
if and only if they are synchronised by the given cospan.  Furthermore, by inspection, the two Forward propagations
coincide, and by symmetry the same inspections shows that the two Backward propagations coincide.
\qed  % (inside enumeration to save space)
\end{enumerate}

The inverted commas around ``pullback'' are to remind us that while the
pullback in $\cat$, along with the lenses constructed on its pullback projections,
might look like a pullback diagram in $\alens$,
it is not necessarily a pullback in that category
(the universally determined mediating functors do not in general have canonical lens
structures on them).

\begin{remark}
\label{rem-spancompn}
{\rm
Given composable representatives of symmetric lenses, that is, spans in $\alens$
which agree on one of their feet (as in $\bV_2 = \bV'_1$ in the proposition),
they can be composed using the ``pullback'' of the cospan exactly as one does for span
composition in $\cat$.  %so there should be a forgetful functor way of saying this...
In more detail:
Imagine that the two triangles in the proposition are the right side of the bow-tie displayed in
Definition~\ref{def-symlens}, and the left side of a similar bowtie in which all the labels have added primes,
then the pullback of the proposition gives a new span of asymmetric lenses with peak $\bT$ say, between $\bS$ and $\bS'$.
As in ordinary span composition these asymmetric lenses can be composed
with the asymmetric lenses $\bS \to \bV_1$ and $\bS' \to \bV'_2$
to yield the composite span of asymmetric lenses
with peak $\bT$ and feet $\bV_1$ and $\bV'_2$.
% Perhaps INSERT double bow tie picture here?  Or not needed?
Here is the picture, labelling the various lenses with their Gets, and in which $H$ and $H'$
are the pullback projections in $\cat$.
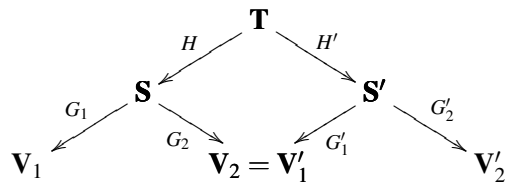
$$
\bfig
\Atriangle(400,250)/->`->`/<400,250>[\bT`\bS`\bS';H`H'` ]
\Vtriangle(400,0)/`->`->/<400,250>[\bS`\bS'`\bV_2=\bV'_1; `G_2`G'_1]
\morphism(400,250)|a|<-400,-250>[\bS`\bV_1;G_1]
\morphism(1200,250)|a|<400,-250>[\bS'`\bV'_2;G'_2]
\efig
$$
Furthermore, using (4) from the proposition, the operational propagations of the composite span agree with
the composite of the operational propagations as described just before the proposition.
%%% This is trying to capture two things I formerly said in a proposition:
%\item This square, which will be referred to as the ``pullback'' of the cospan, can be used to define,
%exactly as ordinary pullbacks do for span composition, a composition of symmetric lenses
%\item The Forward and Backward propagation of the resulting composite symmetric lens correspond
%to the operationally determined propagations of the paragraph before the proposition.
} %\rm
\end{remark}

\section{Symmetric lenses as cospans of asymmetric lenses}
\label{sec-cospan}
%Not always possible. Necessary conditions in JR ZX-paper
%Simplicity and importance for implementation.
%Possibly mention how common it actually is in consultancy work (as long as the `lens' being
%studied hasn't been composed already).
%Show that the symmetric lenses between ABC Frames and XYZ Warehouse, and between XYZ
%Warehouse and XYZ Sales are cospan representable.
%Show that the composite symmetric lense between ABC Frames and XYZ Sales is not cospan
%representable.

In the previous section we saw that every cospan of asymmetric lenses yields, by pullback,
a span of asymmetric lenses, that is, a representative for a symmetric lens.
In fact, the cospan presentation of an asymmetric lens is especially valuable
and is the main way system interoperations are actually built.

To revisit our example from Section~\ref{sec-eg}, we have already seen that
there is an asymmetric lens between ABC Frames and the (possibly imaginary, but frequently
built) system $\bX$.  Furthermore, there is an asymmetric lens from XYZ Warehouse to $\bX$.
The Get of that lens is just a projection from among all the data stored at XYZ Warehouse,
and returns simply the current state of the set $\bX$ in XYZ Warehouse.  The Put starts from a known state of
XYZ Warehouse and a new state of $\bX$, and constructs a new state of XYZ Warehouse by
changing its set $\bX$ to match, leaving the sets $\bY$ and $\bZ$ unchanged,
recalculating the product $\bX \times \bY \times \bZ$, and usually leaving the set of
{\sf Orders} unchanged, but if some elements of $\bX$ have been deleted, and if there
are orders depending on those elements of $\bX$, then those orders are also deleted
(called by database people a ``cascading delete'').  Notice that if instead new elements
of $\bX$ had been inserted, then {\sf Orders} would not change, but the injection into
$\bX \times \bY \times \bZ$ would be adjusted to account for its new larger codomain.

Thus we have a cospan of asymmetric lenses between ABC Frames and XYZ Warehouse.
We could say now, following the previous section, that we ``pullback'' the span to obtain
a representative for a symmetric
lens, thus providing interoperations between ABC Frames' and XYZ Warehouse's systems.  That is indeed
theoretically true.  The resulting system at the peak of the span $\bT$ is sometimes called the
\emph{federated information system} because it is the state space for the system that combines
all of the information held at ABC Frames with all of the information held at XYZ Warehouse,
subject only to ensuring that those two subsystems remain consistent via the same $\bX$ state.
Such symmetric lenses are theoretically important because we can reason with them and prove properties of the
combined system (for example, that certain things remain consistent or that certain operations avoid deadlock or \dots).
But these systems are hardly ever built.  To begin with, ABC Frames and XYZ are separate companies,
and are unlikely to want to, or indeed be able to, break commercial-in-confidence agreements and
share all data that they might hold.  There are commercial, privacy, and cyber security \cite{jscitpmdsd},
reasons, to name just a few, for not building the system $\bT$.

Instead, the system $\bX$ might be built, along with the two asymmetric lenses to it
described above (one from ABC Frames and one from XYZ Warehouse).  Or, alternatively,
the Forward and Backward propagations from such a cospan can be implemented as message
passing and through Applications Programmer Interfaces (APIs) the messages can keep
the two systems of the two companies synchronised.  These options limit the exposure of
each of the companies and their systems to the minimum required
for the system interoperations \cite{jscitpmdsd}, and those system interoperations are in the interests of
the efficiencies of both organisations (after all, we only build such systems if there is a
commercial imperative).

Of course, as noted in Proposition~\ref{prop-lenspb} part 4, the propagations determined
by the cospan through $\bX$ or by the span through the federated information system
are the same.  But the former is a minor piece of engineering work, which can even be
separated into three tasks:  Implementing the small common system $\bX$, and the two
asymmetric lenses from ABC Frames to it (which can be done exclusively by ABC Frames
engineers) and from XYZ Warehouse to it (which can be done exclusively by XYZ Warehouse
engineers).  On the other hand, working with the federated system, either by constructing it
or by simulating propagations through it, is a major piece of work that is generally
hard to partition into secure and independent tasks.

The message of this section is that cospan representations of symmetric lenses
are very much preferred for engineering purposes.

It is worth noting however, that not all symmetric lenses have cospan representations.
The paper \cite{jrcsl} establishes necessary and sufficient conditions for the
existence of cospan representations.  For now, suffice it to say, having a
cospan representation is something that one wants to keep.

And so there is another important point to note:  The composition of symmetric
lenses does not preserve cospan representability.  Two cospan representable symmetric
lenses may compose to give a symmetric lens which is not in itself cospan
representable.

Again, the example from Section~\ref{sec-eg} provides an illustration.

We will not work through the details here, but there is a cospan
representation for the symmetric lens between XYZ Warehouse and XYZ Logistics.
To make the example more realistic we have included in this interoperation
an example of \emph{half-duplex interoperation} (see \cite{djhdi}).  In short,
XYZ Logistics is not permitted to change the state of the catalogue {\sf XYZ} --- it
is read-only.  There is also an opportunity here to introduce a
non-trivial amendment lens (see \cite{dml}) between XYZ Warehouse and XYZ Logistics,
but to keep things simple let's assume that the company XYZ enters, processes and fills single orders
at a time (otherwise orders in XYZ Logistics might have to be reversed (amended) by XYZ Warehouse
if the monic constraints in XYZ Warehouse were violated).

The cospan of symmetric lenses between XYZ Warehouse and XYZ Logistics determines
by ``pullback'' a representative for a symmetric lens.  The two symmetric lenses
(between ABC Frames and XYZ Warehouse, and between XYZ Warehouse and XYZ Logistics)
can be composed, either by ``pulling back'' (creating an even larger federated system
$\bT''$), or by composing the propagations, and it may be that for whatever reason
the composite symmetric lens is the subject of our interest.  But note well:  The
composite symmetric lens is not cospan representable.  Presented merely with the
composite symmetric lens (and so, no information about how XYZ Warehouse mediates
the information between ABC Frames and XYZ Logistics) there is no simple shared
data that the two organisations can synchronise upon.  The super federated system
could be used in theory to build interoperations, but the information about the
engineering appropriate cospans is gone.

Perhaps it would be better if the example from Section~\ref{sec-eg} were treated
as a multilens, since then all three organisations, and their interactions, could
be captured in a single mathematical entity.

\section{Multilenses}
%Mention `amendment', but note that it is not needed in this paper, and that `ordinary' multilenses are a simple
%special case of the `amendment' multilenses presented by Diskin and treated in our multilens paper.
%Cite Diskin and JR Multilens paper.

\begin{definition}For $n$ a positive integer, an \emph{$n$-lens} consists of %a category $A_1$ and
$n$ asymmetric lenses $f_i$ with common domain $\bS$, such that $f_i : \bS \to A_i$.
\end{definition}

%Note that if $n>0$ then $f_1$ has the distinguished category $A_1$ as its lens codomain.  A 0-lens is simply a category $A_1$.
A 1-lens is an asymmetric lens $f_1 : \bS \to A_1$.  A 2-lens is a representative for a symmetric lens ---
a span in $\alens$ as in Definition~\ref{def-symlens}.

For an $n$-lens $L = (f_i : \bS \to A_i)$, since it is in general an $n$-wide span in $\alens$, we adopt, and adapt,
the terminology usually used for parts of wide spans (including ordinary spans).   Thus the $f_i$ are called the
\emph{legs} of $L$, and the $A_i$ are called the \emph{feet} of $L$.  The category $\bS$ is called the \emph{peak} of $L$.
%When $n>0$ w
We call $f_1$ the
\emph{leftmost} leg of $L$ and $f_n$ the \emph{rightmost} leg of $L$, and $A_1$ the \emph{leftmost} foot
of $L$ and $A_n$ the rightmost foot of $L$.  Of course, for $n=1$ the leftmost and rightmost legs of $L$
coincide, and are both $f_1$, and likewise $A_1$ is both the leftmost and the rightmost foot of $L$.
%When $n=0$ there are no legs at all, but $A_1$ is still called a foot, and may be referred to as both
%the leftmost and the rightmost foot.
When $n \le 2$, the ``most'' of leftmost or rightmost is superfluous in
normal usage, and it is common to say just ``the left leg'' or ``the right foot'' etc, and even when $n>2$,
if there is little chance of confusion, we may still say ``left'' for ``leftmost'' and ``right'' for ``rightmost''.

For $n>1$ an $n$-lens is a wide span in the category $\alens$. % with distinguished left foot.
Relating this to previous work, an $n$-lens for $n>1$ is
a multiary lens \cite{jrmml} in which every asymmetric amendment lens is closed,
that is, all amendments are trivial.
Such $n$-lenses form the ``special case'' (wide spans of d-lenses) % (as opposed to wide spans of asymmetric
%amendment lenses), %(of which d-lenses are the closed case),
referred to in the final paragraph of \cite{jrmml}.
Thus, for $n>1$, $n$-lenses are a specialisation of multiary lenses --- the special case in which all
amendments are trivial.  It may be worth emphasising that this ``special case'' is what the authors see
as the main case.  There are occasionally circumstances in which non-trivial amendments are useful, and
the paper \cite{jrmml} dealt with nontrivial amendments to have the broadest possible generality and to link
directly with the extant work of Diskin et al \cite{dml}, but in this paper we restrict our attention
to multilenses: wide spans of asymmetric lenses without amendments.

\begin{definition}
A \emph{multilens} $L$ is an $n$-lens for some (positive integer) $n$.
If $n>1$ the multilens $L$ is said to be \emph{non-trivial}.
\end{definition}

\section{The fusion of multilenses}
%Define fusion --- an $m+n-1$ multilens instead.  State and possibly prove any easy enough properties.
%Note that pairwise, and even $n$-wise interoperations among legs of a multispan may be (wide-) cospan representable. -- NOT DONE
%State and perhaps `prove' that fusion preserves wide-cospan representability.  (I don't think there's much proof required.) -- NOT DONE

The multiary lenses of \cite{jrmml} compose, as shown there, with a multicategory \cite{lhohc} structure.  In the
terminology of this paper, using the composite defined in \cite{jrmml}, an $m$-lens and an $n$-lens compose
to give an $(m+n-2)$-lens (think for example of ``plugging'' the left leg of one lens into the right
leg of the other with those two legs ``disappearing'').  That composition is a generalisation of the
usual composition of symmetric lenses, or indeed of spans or relations --- a $2$-lens composes with a
$2$-lens if the leftmost foot of one equals the rightmost foot of the other,
and the result is a $(2+2-2)$-lens with peak a pullback calculated over the common foot
(see Remark~\ref{rem-spancompn}).  Notice that in this familiar composition
the common foot and the two legs to that foot all disappear
(hence the subtraction of two in the count of legs).

The simple, but important, change in this paper is the introduction of a new composition called
fusion which retains the foot that has been composed over. %in order to preserve information about the factorisation.

\begin{definition}
\label{def-fusion}
Suppose $L = (f_i : \bS \to A_i)$ is an $m$-lens and
$L' = (f'_i : \bS \to A'_i)$ is an $n$-lens
with the rightmost foot of $L$ being equal to the leftmost foot of $L'$.
Then the \emph{fusion} of $L$ and $L'$, denoted here simply by the juxtaposition $LL'$,
is the lens $LL' = (g_i : \bS \to B_i)$ given as follows:
$LL'$ is an $(m+n-1)$-lens with feet
  $B_i = A_i$ for $i \le m$, and
  $B_i = A'_{i-m+1}$ for $i \ge m$.
Let $\bT$ be the pullback of $f_m$ along $g_1$ with projections
$H$ and $H'$, then
$g_i = f_i H$ if $i \le m$, and
$g_i = f'_{i-m+1} H'$ if $i \ge m$.

%If $m$ (or $n$) is zero then $g_i = f'_i$ (or $g_i = f_i$)  and $LL'$ is an $n$- (or $m$-) lens (respectively).
%If both $m$ and $n$ are zero then $LL'$ is a $0$-lens with distinguished foot $B_1 = A_1 = A'_1$.
\end{definition}

\begin{remark}
{\rm
We record here a few basic results about the fusion operation.
\begin{enumerate}
\item Well-definedness:  The use of both $i \le m$ and $i \ge m$ in the definition is deliberate, and is intended to
reinforce the sense of fusion.  %Writing in the general case where $n$ and $m$ non-zero,
If $i=m$ then
$B_i = A_i = A'_{i-m+1}$ by assumption, and $g_i = f_i H = f'_{i-m+1} H'$ by Proposition~\ref{prop-lenspb} part 2,
so the fusion is well-defined.
\item Identities:  Identity $1$-lenses are, up to equivalence, left and right identities for fusion.
\item Associativity:  Up to span isomorphism in $\cat$, the fusion operation is associative.  The equivalence
relation presented in \cite{jrjot} (and mostly avoided here) is coarser
than span isomorphism in $\cat$, and is a congruence for fusion (and for the composition of \cite{jrmml}), so
fusion is also associative for equivalence classes of multilenses.  (The calculations are tedious, but routine, and
follow the path traced in \cite{jrjot}, so they have been suppressed here.)
\end{enumerate}
} %\rm
\end{remark}

We would like to emphasise that fusion is a minor change from multilens composition.  Non-trivial multilenses are
fusable if and only if they are composable ---
fusion simply keeps the foot that one composes over along with the single (by Proposition~\ref{prop-lenspb} part 2)
leg to that foot.  It is still an operation which combines composable multilenses to get multilenses.
But fusion often feels like a significant change for people who are used to composing symmetric lenses
because the fusion of two symmetric lenses is not a symmetric lens but rather a 3-multilens.
This difference is exactly what we need for our applications.  We will return to this in Section~\ref{sec-otherfusions}
where we illustrate a few other well-known fusion operators for comparison purposes and to set readers'
minds at ease.

\section{Sometimes lens fusions are better than lens compositions}
\label{sec-fusionsbetter}
%Revisit example to see these things illustrated.
%Basic examples first:

We turn now to some basic examples of fusion, and then revisit the example of Section~\ref{sec-eg}.

What happens when we fuse $1$-lenses, recalling that $1$-lenses are themselves
simply asymmetric lenses?

If $L$ and $L'$ are $1$-lenses, then $LL'$ is also a $(1+1-1=1)$-lens, so fusion is an operation on asymmetric
lenses.  But it is not the usual composition of asymmetric lenses because fusable 1-lenses have common
codomains.  That is, they form a cospan of asymmetric lenses $L: \bS \to A_1 = A'_1 \toleft \bS' : L'$.
So, how do we fuse a cospan?
Definition~\ref{def-fusion} tells us that we pull the two asymmetric lenses back along each other,
and the resulting asymmetric lens is the diagonal of the ``pullback'' square $\bT \to A_1$.  This is sometimes
known as the \emph{consistency lens}.  In database terms, if the trough of the cospan $A_1$ is the system of states
of common data, then the peak of the ``pullback'' is, up to isomorphism, the category whose objects are
consistent pairs of states of the systems $\bS$ and $\bS'$, consistent in as much as they share the same
common data state, and whose arrows are pairs of transitions,
one from $\bS$ and one from $\bS'$ which are consistent in as much as they involve the same transition
in $A_1$ for the shared data.  The $1$-lens $LL'$ is an asymmetric lens.  The Get, the functor part of the
diagonal, tells us how the shared data changes when a $\bT$ transition takes place, and the Put tells us
how to change the consistent states in $\bT$ when the shared data is changed.

But there is yet another way that we might fuse the asymmetric lenses $L$ and $L'$.
It is well-known that an asymmetric lens can be represented as a symmetric lens in two ways:
For the asymmetric lens $L$, form a span of asymmetric lenses (Definition~\ref{def-symlens})
by pairing $L$ with the identity on $\bS$ on either the left of the right.
If we do that on the left for $L$, and on the right for $L'$ (using of course the identity on $\bS'$)
we obtain two $2$-lenses which we know are symmetric forms of the asymmetric
lenses $L$ and $L'$, and these two $2$-lenses are again fusable. This time Definition~\ref{def-fusion}
tells us that the resulting lens with be a $(2+2-1=3)$-lens:  It is the three legged ``pullback'' cone over the
cospan $\bS \to A_1 = A'_1 \toleft \bS'$ consisting of the consistency lens in the middle, and the two
``pullback'' projection lenses $H : \bT \to \bS$ and $H' : \bT \to \bS'$ on the left and right.

The last paragraph describes a special case of the fusion of two symmetric lenses:  Two $2$-lenses
fuse to form a $3$-lens.  The three legs are again the consistency lens in the middle, and the outer two
legs are, together, the usual symmetric lens composite.  The fusion ``remembers'' the foot $A_1$ that
the symmetric lenses have been composed over, and its relationship to the peak $\bT$ via the
consistency lens.  This is a small, but important difference, as we will soon see.  The fusion
remembers the way the composed up symmetric lens factors into two symmetric lenses.  And why is this
important?  It is because, as noted in Section~\ref{sec-cospan}, the fact that there might be a cospan
representation of each symmetric lens is important in engineering, but the composed symmetric lens
might have no cospan representation so the factorisation is vital for actually building system
interoperations.  Let's look again at the example presented in Section~\ref{sec-eg} and further
developed in Sections~\ref{sec-asymlens} and~\ref{sec-cospan}.

Recall that there are symmetric lenses, $2$-lenses, between ABC Frames and XYZ Warehouse, and
between XYZ Warehouse and XYZ Logistics, and that those two symmetric lenses are cospan representable,
but that the composite $2$-lens between ABC Frames and XYZ Logistics is not cospan representable.

If we are concerned, as was hypothesised in Section~\ref{sec-cospan}, with the interoperations between
ABC Frames and XYZ Logistics, we might have composed the two symmetric lenses and lost the information
that the composite factors into two symmetric lenses that are cospan representable.  But alternatively,
we could have fused the two symmetric lenses, the two $2$-lenses, to obtain a $3$-lens.  That $3$-lens
does indeed contain all the information required to study and prove properties about the interactions
between ABC Frames and XYZ Logistics.  But it also includes the base information to see how to factorise
that interaction through XYZ Warehouse, and through that factorisation and the cospan representations,
the interoperations are then easy to engineer as well.

Sometimes, lens fusions are better than lens compositions.

%Fusion of asymm lenses as symmetric lenses $2+2-1=3$.
%Asym lens cospan arrows seem to have disappeared, but are still contained in the span to the middle.
%Compare our example -- composite overall symm lens has lost the middle and has no way of reconstructing it,
%or even just the normal composite of a cospan of asym lenses as sym lenses --- the factorisation is lost.

Thus, the more legs the merrier:  In our example in fact, if we have the cospans the best fusion is five legged!
See below.  But first, a little more about fusions.

\section{Other fusions}
\label{sec-otherfusions}
%Provided there is space:  Describe other similar operations including

At first sight, some people find the fusion operator confronting because it works like a composition
(at least, one can only fuse composable lenses), but $n$-lenses, and especially concerningly the familiar
$2$-lenses normally called symmetric lenses, aren't closed under fusion.  In general, fusion takes things of certain
types and produces things of different types.

So it seems worthwhile to point out that operators like fusion are common, and have a long history in
mathematics and software engineering.  To offer just a few examples:
\begin{enumerate}
\item Path categories:  Paths in a topological space $X$ are normally defined to be continuous functions
from the unit interval $[0,1]$ to $X$.  Two paths with common end and start points can be composed by reparameterization
to get a new path --- the first path is traversed ``twice as fast'' so that the second path can be traversed in
the second half of the unit interval.  Probably this is all familiar.  Of course, such composition is not associative
because three paths will be traversed in either the first two quarters and the last half of the unit interval or
the first half and the last two quarters of the unit interval depending on which composition is taken first.  So,
we introduce an equivalence relation, homotopy, which allows reparameterized paths to be treated as equivalent.

Alternatively, one might approach the problem in the style of fusion.  Let paths be continuous functions from
an interval of length $n$, say for non-negative integers $n$, and paths are fused (but usually still called
``composed'') by having a path from $[0,m]$ and a path from $[0,n]$ form a path from $[0,m+m]$.  This fusion
is associative and has paths of length zero, paths from $[0,0]$, also known as points, as identities.
Such paths form a category, the Moore Path Category, with no need for any equivalence relation.
\item Free monoid construction (including list concatenation):   The discrete form of the example just given is
familiar to computer scientists.  Lists, or indeed arrays, might be seen as functions from sets $[n]$ of
$n$ elements into a data type $X$.  The concatenation of lists of length $m$ and length $n$ gives a list of length
$m+n$ by fusing the functions.  In mathematics of course, this is the construction of the free monoid on a set $X$
via words in the alphabet $X$.  Again associativity is immediate, and identity comes from the empty domain $[0]$.
\item Other free constructions from $F1$:  In both the fusion examples just given, the collection of domains
has the form of a coalgebra family \cite{jwaoalfflt} and can be calculated from the free algebra on the terminal.
There are more examples of the same kind.
%\item Parametric right adjoints more generally (check!)
\item Composition of (lax) natural transformations:  To offer an example of a rather different kind, natural
transformations are, in chain complex terms, degree 1 maps.  Given a natural transformation
$\eta : F \to G : \cal A \to \cal B$ an \emph{object} of $\cal A$ is sent to an \emph{arrow} $\eta _A$ of $\cal B$, and an
\emph{arrow} of $\cal A$ is sent to a commutative \emph{square} in $\cal B$, or if $\cal B$ is a $2$-category and $\eta$ a
lax natural transformation, to a \emph{$2$-cell} of $\cal B$, etc.  The horizontal composite of natural transformations
is normally said to yield a natural transformation because a commuting square is indeed a $1$-cell.  But a
very natural composition, in analogy with the examples above, one might say a fusion, of lax natural transformations,
yields a modification (a degree $2$ map), amd indeed the fusion of three lax natural transformations yields a
perturbation (a degree $3$ map).  And more generally degree $m$ and degree $n$ maps in higher categories can
fuse to give degree $m+n$ maps.  An extensive study of algebras with these kinds of compositions was undertaken by the
Dutch mathematician Sjoerd Crans under the name \emph{teisi} (singular \emph{tas}).
\end{enumerate}

We end this section with a note about counting.  Recall that the fusion of an $m$-lens and an $n$-lens is an
$(m+n-1)$-lens.  But the examples just given ``fuse'' $m$-structures and $n$-structures to get $(m+n)$-structures.
In fact there is no substantive difference, and the apparent difference arises just from how we count.  The natural
way of counting, and hence of labelling, multilenses is by counting the legs.  But the examples above are counted and
labelled instead by the equivalent of the \emph{spaces between} the legs (the interval $[0,n]$ for example is of
length $n$ because there are $n$ unit intervals in the spaces between the ``feet'' $\{0,1, \dots n\})$.  The formulas coincide
exactly if we use the same counting paradigm in each case.

A further example, and one in which the natural count again corresponds to the way we count legs
in $n$-lenses, is the natural join operation for databases (this example was suggested to us by an anonymous referee).
To take it in its simple form, two tables in a database consisting of say $m$ and $n$ columns respectively,
are joined on some common data (often a key attribute) to yield a table with $m+n-1$ columns.

%\section{Related work}
%Describe more clearly history.  Cite Zinovy.  Cite K\"onig. Cite micro-services German guy.
%
%In the long version of this paper we include historical remarks on related work and
%more recent applications of the ideas outlined here.

\section{Conclusion and further work}
%In short, sometimes multilenses are better than composed up symmetric lenses.

We have seen that the \emph{fusion} operation is sometimes better than the normal composition of, for example,
symmetric lenses, because it preserves information about what structures have been composed over, and this factorisation
information may be very valuable in knowing where to find cospan representations. And those cospan representations may
be very useful in the implementation of systems interoperations.

But, the reader might ask, if we have those cospan representations, why would we even do fusion?  We don't want
to lose the cospan representations if they're so useful for implementations, and the cospan representations are
\emph{among} the feet of the spans.  The fusion preserves legs and feet, but nothing between the feet.  The results of
fusions are wide spans: A peak, some legs, and some (bare) feet.

In further work we have shown %via a study of the properties of the category $\alens$ of asymmetric lenses
that this loss of the cospans need not be an issue at all.

If we have the cospan representations, then we begin with a ``zig-zag`` of asymmetric lenses among the feet.
As discussed in Section~\ref{sec-fusionsbetter} we can treat each of those asymmetric lenses as a symmetric lens,
a $2$-lens, (by pairing it, on the \emph{outside} of the cospan, with identities) and then fuse the $2$-lenses
so that each cospan becomes a $3$-lens.  We can even do all this in one go --- and this is a general result ---
and the fusion will be a wide span of asymmetric lenses canonically built on the limit cone of the zig-zag as calculated
in $\cat$.

In the case of our supply chain example, the previous paragraph means that we end up with a $5$-lens. The five legs
are asymmetric lenses with codomains the three business systems (ABC Frames, XYZ Warehouse and XYZ Logistics) and the
two common data subsystems ({\sf X} and ${\sf Orders \to XYX}$).  As noted above, the cospans among the feet are indeed
gone, but their objects remain, and, remarkably, the cospans can be \emph{uniquely} recovered whenever required.
This follows from a particularly nice orthogonal factorisation system on $\alens$ using image factorisations.  And
image factorisations are much simpler for lenses than for arbitrary functors, and are frequently used, often
tacitly, by engineers who cut lenses down to their image factorisations routinely.  The unique fillers for the
orthogonal factorisation system restore the asymmetric lenses making up the cospans whenever they are required.
The fusion contains not just the factorisation information required to search for the cospan representations,
but in fact all the information needed to fully determine the cospan representations, and so all the information
needed for an effective implementation.

%Include the image orthogonal factorisation system on $\cat$, and how that gives effectively UFAL (unique factorisation of
%asymmetric lenses), so that a wide span with merrily many legs determines uniquely the cospans between its feet
%that are needed for actually building the interoperations in practice.

%Also, remark upon the multiple fusion of multiple cospan representable lenses being the limit cone
%over the cospan diagram.

\section{Acknowledgements}
The authors gratefully acknowledge the support of the Australian Research Council,
the Centre of Australian Category Theory and Mount Allison University.  We have
benefited from the perspicacious comments of anonymous referees and from
useful discussions with Bryce Clarke and Angus Johnson, and we thank
them all for sharing their insights.

\end{document}